\documentclass{article}

\usepackage{amssymb}
\usepackage{amsmath}
\usepackage{amsthm}
\usepackage{mathrsfs}
\usepackage{cite}
\usepackage{bbm}

\usepackage{pgf,tikz}
\usepackage{tikz-cd}

\newtheorem*{theorem*}{Theorem}
\newtheorem*{lem*}{Lemma}

\begin{document} 

\begin{center}
\begin{Large} {\bf The complex plank problem, revisited} \\[0.4cm] \end{Large}

\begin{large} Oscar Ortega-Moreno \end{large}
\end{center}

\vspace{-0.9cm}

\begin{quote}
\footnotesize{ \vskip 1cm \noindent {\bf Abstract.}
Ball's complex plank theorem states that if $v_1,\dots,v_n$ are unit vectors in $\mathbb{C}^d$, and $t_1,\dots,t_n$, non-negative numbers satisfying $\sum_{k=1}^nt_k^2 = 1,$ then there exists a unit vector $v$ in $\mathbb{C}^d$ for which $|\langle v_k,v \rangle | \geq t_k$ for every $k$. Here we present a streamlined version of Ball's original proof.  
}
\end{quote}

\vspace{0.4cm}

\setlength{\parindent}{0cm}
\noindent

A \emph{plank} in a vector space is the region between two parallel hyperplanes. Ball's plank theorem \textbf{\cite{Ball1991}} states the following:

\begin{theorem*}
If the unit ball of a Banach space $X$ is covered by a (countable) collection of planks in $X$, then the sum of the half-widths of these planks is at least $1$.
\end{theorem*} 

For Hilbert spaces, the statement was proven almost 40 years earlier by Bang \textbf{\cite{Bang1951}} in his remarkable solution to the famous plank problem of Tarski from the 1930s.\\
As it stands, the plank theorem is sharp: any ball of a Banach space may be covered by a collection of non-overlapping parallel planks whose half-widths add up to 1. However, if one considers only planks which are symmetric about the origin, the condition of the theorem may be improved depending on the geometry of the space $X$. For example, for complex Hilbert spaces, the complex plank theorem \textbf{\cite{Ball2001}} (stated in the abstract) provides the optimal condition: the sum of the squares of the half-widths of the planks should be at least $1$. \\For real Hilbert spaces, perhaps the most natural way to state the optimal condition is in the form of Fejes T\'oth's zone conjecture. A zone of spherical width $w$ is the intersection of the sphere $\mathbb{S}^{d-1}$ and an origin symmetric plank of width $\sin(w/2)$. In 1973, Fejes Tóth \textbf{\cite{Fejes1973}} conjectured that if $n$ zones of equal angular width cover the sphere then their angular width is at least $\pi/n$. This conjecture was recently resolved by Jiang and Polyanskii \textbf{\cite{JP2017}} whose proof is based on the ideas of Bang\textbf{\cite{Bang1951}} and Goodman and Goodman \textbf{\cite{GG1945}}. In \textbf{\cite{OM2021}} the author of the present note found a completely different proof inspired by Ball's solution to the complex plank problem.\\
Recently, Zhao \textbf{\cite{Zhao2021}} simplified the author's proof of the zone conjecture by eliminating the reformulation of the problem in terms of Gram matrices and inverse eigenvectors.\\
In this note, we simplify Ball's original proof of the complex plank problem in a similar way, thereby obtaining a direct proof of it. We start with the case when all the $t_k$ are equal to $\frac{1}{\sqrt{n}}$.

\begin{lem*}
Let $v_1,\dots,v_n$ be unit vectors in $\mathbb{C}^d$. If $u$ maximizes $\prod_{k = 1}^n|\left\langle v_k, u \right\rangle|$ among unit vectors, then
\begin{equation}\label{eq1}
	u = \frac{1}{n}\sum_{k = 1}^n\frac{1}{\left\langle v_k , u \right\rangle}v_k.
\end{equation} 
\end{lem*}
 
The proof of this lemma is an immediate application of Lagrange multipliers. By compactness of the unit sphere in $\mathbb{C}^d$ it is enough to show the following theorem.

\begin{theorem*}
Let $v_1,\dots, v_n$ be unit vectors in $\mathbb{C}^d$. If $u$ maximizes $\prod_{k=1}^n|\langle v_k ,u \rangle|$ among unit vectors, then $|\langle v_k ,u \rangle|\geq \frac{1}{\sqrt{n}}$ for every $k$. 
\end{theorem*}
\noindent {\it Proof.} Suppose for a contradiction that $|\langle v_1 ,u \rangle|<\frac{1}{\sqrt{n}}$ (note that $|\langle v_k ,u \rangle|>0$ for every $k$ due to the choice of $u$). We will examine the values of the product on vectors of the form 
\[
u_z = z v_1 + u
\]
for different complex numbers $z$ such that $u_z$ belongs to the complex unit sphere. Note that $\|u_z\|^2 = |z|^2 +2 \mathfrak{R}(z \langle v_1 , u \rangle) + 1.$ So $u_z$ belongs to the unit sphere if and only if $|z + \langle u ,v_1 \rangle| =|\langle v_1 , u \rangle|$, i.e. $z$ belongs to the circle of radius $|\langle v_1 , u \rangle|$ centered at $-\langle u ,v_1 \rangle$. Let us call this circle $C$. The first factor of the product is also constant on $C$ with value  $|\langle v_1 ,u \rangle|$, since
$
|\langle u_z ,v_1 \rangle| = | z + \langle u ,v_1 \rangle| = |\langle v_1 ,u \rangle|.
$
Let $p$ be the complex polynomial defined by 
\[p(z) = \prod_{k=2}^n\frac{\langle u_z, v_k  \rangle}{\langle u , v_k \rangle} = \prod_{k=2}^n\frac{\langle z v_1 + u, v_k  \rangle}{\langle u ,v_k \rangle}. 
\]
Clearly, the maximum of $|p(z)|$ on the circle $C$ is $p(0) =  1$ since $u_0 = u$. By (\ref{eq1}),
\[p'(0) = \sum_{k = 2}^n\frac{\langle v_1 ,v_k \rangle}{\langle u ,v_k \rangle} = \left\langle  v_1,nu-\frac{v_1}{\langle v_1,u \rangle} \right\rangle = \langle v_1 , u \rangle \left(n - \frac{1}{|\langle v_1 , u \rangle|^2}\right).
\]
Note that $n - \frac{1}{|\langle v_1 ,u \rangle|^2 } < 0.$ Therefore, by moving along the direction $-\langle u,v_1 \rangle$ from $0$ into the open disc enclosed by the circle $C$, we find a complex number $z$ in the open disc with $|p(z)| > 1$. By the maximum modulus principle, there exists a complex number $z$ in $C$ with $|p(z)|>1$, contradicting the choice of $u$. 
\hfill $\square$

\vspace{0.3cm} 

For the general case, we choose $u$ in a slightly different way by maximizing the product $\prod_{k = 1}^n|\left\langle v_k,v \right\rangle|^{t_k^2}$ among unit vectors $v$.

\begin{lem*}
Let $v_1,\dots,v_n$ be unit vectors in $\mathbb{C}^d$ and $t_1,\dots,t_n$, positive numbers satisfying $\sum_{k=1}^nt_k^2 = 1$. If $u$ maximizes $\prod_{k = 1}^n|\left\langle v_k,u \right\rangle|^{t_k^2}$ among unit vectors, then
\begin{equation*}
	u = \sum_{k = 1}^n\frac{t_k^2}{\left\langle v_k ,u \right\rangle}v_k.
\end{equation*} 
\end{lem*}
The rest of the proof is the same as the one for equal width planks, except that one has to analyze the function 
\[
p(z) = \prod_{k=2}^n\left|\frac{\langle v_k ,u_z \rangle}{\langle v_k ,u \rangle}\right|^{t_k^2}
\]
on $C$. Since the logarithm of $p$ is subharmonic, the maximum principle is still valid. 

\vspace{1cm}
\noindent {{\bf Acknowledgments} The author was supported by the Austrian Science Fund (FWF), Project number: P31448-N35.

\begin{small}
\noindent Oscar Ortega-Moreno\\
Vienna University of Technology\\
oscarortem@gmail.com
\end{small}

\end{document}